\documentclass[11pt]{amsart}
\usepackage{amsmath,amssymb,amsthm,mathtools}
\usepackage{hyperref,enumitem,xcolor, cite}
\hypersetup{hidelinks}
\usepackage[margin=1in]{geometry} 

\newtheorem*{theorem}{Theorem}

\title[Isometric copies of $\ell_\infty^n$ in Lipschitz-free spaces over finite metric spaces]{Isometric copies of $\ell_\infty^n$ in Lipschitz-free spaces\\ over finite metric spaces}

\author{Rainis Haller}
\address{Institute of Mathematics and Statistics, University of Tartu, Estonia}
\email{rainis.haller@ut.ee}

\thanks{This work was supported by the Estonian Research Council grant PRG1901.}

\begin{document}

\begin{abstract}
For every $n\in\mathbb N$, we construct a finite metric subspace $M_n$ of $\ell_\infty^n$ such that the Lipschitz-free space $\mathcal F(M_n)$ contains a linear isometric copy of $\ell_\infty^n$. This answers a question posed by Khan, Mim, and Ostrovskii, who obtained examples in dimensions three and four. 
As a consequence, the finite-dimensional real Banach spaces that admit a linear isometric embedding into a Lipschitz-free space over a finite metric space are precisely the polyhedral spaces.
\end{abstract}

\maketitle

\section{Introduction}

Let $M$ be a pointed metric space, and let $\mathcal F(M)$ denote its Lipschitz-free space. 
For finite $M$, the space $\mathcal F(M)$ is polyhedral, and so are all its subspaces. We ask whether the converse holds: does every finite-dimensional polyhedral Banach space embed linearly and isometrically into $\mathcal F(M)$ for some finite metric space $M$? It suffices to answer this question for $\ell_\infty^n$, since every finite-dimensional polyhedral Banach space embeds linearly and isometrically into $\ell_\infty^n$ for some $n$.
In \cite{MR4312039}, Khan, Mim, and Ostrovskii studied isometric copies of $\ell_1^n$ and $\ell_\infty^n$ in Lipschitz-free spaces over finite metric spaces. In the case of $\ell_\infty^n$, they constructed finite metric spaces $M_3$ and $M_4$ such that $\mathcal F(M_3)$ and $\mathcal F(M_4)$ contain linear isometric copies of $\ell_\infty^3$ and $\ell_\infty^4$, respectively. They asked whether such a finite metric space exists for $n\geq 5$. 

Alexander, Fradelizi, Garcia-Lirola, and Zvavitch obtained a related low-dimensional result. They proved that $\mathcal F(M)$ itself can be isometric to $\ell_\infty^n$ only for $n\leq 3$; when $n=3$, this occurs precisely when the canonical graph associated with $M$ is a four-cycle with equal edge weights \cite[Corollary~3.8]{MR4181163}. Dilworth, Kutzarova, and Ostrovskii later constructed simple six- and eight-point metric spaces whose Lipschitz-free spaces contain isometric copies of $\ell_\infty^3$ and $\ell_\infty^4$, respectively; see \cite[Section~8]{MR4301142}. Ostrovska and Ostrovskii subsequently obtained several obstructions to the presence of isometric copies of $\ell_\infty^n$ in Lipschitz-free spaces over finite metric spaces; see \cite{MR4455126}.

The problem was brought to the author's attention during the 2026 Summer School \emph{Topics in Banach Space Theory} in Castro Urdiales. During the meeting, the author found similar explicit constructions in dimensions three to six. The construction below arose from an attempt to understand and generalise the common idea behind these examples. The main result is the following.

\begin{theorem}\label{thm}
For every $n\in\mathbb N$, there exists a finite subset $M_n$ of $\mathbb Z^n$ containing $0$, equipped with the metric inherited from $\ell_\infty^n$, such that the Lipschitz-free space $\mathcal F(M_n)$ contains a linear isometric copy of $\ell_\infty^n$. In fact, there exist elements $u_1,\dotsc,u_n\in\mathcal F(M_n)$ such that
\[
\big\| \sum_{i=1}^n t_i u_i \big\| = \max_{1\leq i\leq n}|t_i|\quad \text{
for all $t_1,\dotsc,t_n\in\mathbb R$}.\tag{1}\label{eq:1}
\]
\end{theorem}




Throughout, all Banach spaces are real, and $\delta_x$ denotes the canonical element of $\mathcal F(M)$ corresponding to $x\in M$.

\section{Proof of the theorem}

Fix $n\in\mathbb N$. For $j=1,\dotsc,n$, denote by $\pi_j$ the $j$-th coordinate functional on $\mathbb R^n$. 
Suppose that $M$ is a finite metric subspace of $\ell_\infty^n$ containing $0$, and that $u_1,\dotsc,u_n\in\mathcal F(M)$ satisfy $\langle u_i,\pi_j\rangle=\delta_{ij}$ for all $i,j\in\{1,\dotsc,n\}$ and, for all $\theta_1,\dotsc,\theta_n\in\{-1,1\}$, 
$\|\sum_{i=1}^n \theta_i u_i\|\leq 1$. Then \eqref{eq:1} holds. Indeed, since every $\pi_j$ is $1$-Lipschitz,
\[
\big\|\sum_{i=1}^nt_i u_i\big\|\geq\max_{1\leq j\leq n}\big|\big\langle \sum_{i=1}^n t_i u_i,\pi_j\big\rangle\big|=\max_{1\leq j\leq n}|t_j|.
\]
The reverse inequality follows from the fact that $[-1,1]^n$ is the convex hull of $\{-1,1\}^n$, together with convexity and homogeneity of the norm. 

Let $E=\{e_1,\dotsc,e_n\}\cup\{-1,1\}^n$. We seek non-negative compactly supported functions $\varphi_y\colon\mathbb R^n\to\mathbb R$, $y\in E$, such that
$\sum_{a\in\mathbb Z^n} \varphi_{y}(a)=1$. We shall set
\[
u_i=\sum_{a\in\mathbb Z^n} \varphi_{e_i}(a)(\delta_{a+e_i}-\delta_a)
\]
and require that, for every $\theta=(\theta_1,\dotsc,\theta_n)\in \{-1,1\}^n$, 
\[
\sum_{i=1}^n \theta_i u_i=\sum_{a\in\mathbb Z^n} \varphi_{\theta}(a)(\delta_{a+\theta}-\delta_a).
\]
Comparing the coefficients of $\delta_a$ in $\sum_{i=1}^n\theta_i u_i$ gives
$
\varphi_{\theta}(a-\theta)-\varphi_{\theta}(a)=\sum_{i=1}^n\theta_i \big(\varphi_{e_i}(a-e_i)-\varphi_{e_i}(a)\big).
$
For $y\in\mathbb R^n$, define the difference operator $\Delta_y f(x)=f(x)-f(x-y)$. Then the condition takes the form
\[
\Delta_\theta \varphi_{\theta}=\sum_{i=1}^n\theta_i\Delta_{e_i} \varphi_{e_i}.
\]
For $y\in\mathbb R^n$, write $D_yf(x)=\left.\frac{d}{ds}f(x+sy)\right|_{s=0}$ whenever this derivative exists. Finite differences are not linear in the direction: in general, $\Delta_{y+z}f\neq \Delta_y f+\Delta_z f$. By contrast, for a $C^1$-function, directional differentiation is linear in the direction. 
This suggests constructing a $C^1$-function $\psi$ such that, for every $y\in E$, $D_y \psi =\Delta_y \varphi_{y}$. The required identity then follows from $\theta=\sum_{i=1}^n\theta_i e_i$.

For $y\in \mathbb R^n$ and $f\in C_c(\mathbb R^n)$, define
$(A_yf)(x)=\int_0^1 f(x-ty)\, dt$. Equivalently, $A_y f=f*\mu_y$, where $\mu_y$ is the uniform probability measure on the line segment $[0,y]$. Then $D_y(A_y f)=\Delta_y f$. Indeed, a change of variables gives $(A_yf)(x+sy)=\int_{-s}^{1-s} f(x-ty)\,dt$. The operators $A_y$ are positive, map compactly supported functions to compactly supported functions, and commute. If $\sum_{a\in\mathbb Z^n} f(x-a)=1$ for every $x\in \mathbb R^n$, then
\[
\sum_{a\in\mathbb Z^n} (A_yf)(x-a)=\int_0^1\sum_{a\in\mathbb Z^n} f(x-ty-a)\, dt=1.
\]

Let $h(t)=\max\{1-|t|,0\}$ and put $f_0(x)=h(x_1)\cdot\dotsc\cdot h(x_n)$, $x=(x_1,\dotsc,x_n)\in\mathbb R^n$. Then $\sum_{a\in\mathbb Z}h(t-a)=1$, and hence 
\[\sum_{a\in\mathbb Z^n} f_0(x-a)=\sum_{a_1,\dotsc,a_n\in\mathbb Z} h(x_1-a_1)\cdot\dotsc\cdot h(x_n-a_n)=\prod_{i=1}^n \sum_{a_i\in\mathbb Z} h(x_i-a_i)=\prod_{i=1}^n1=1.\]
Let 
\[
\psi=\big(\prod_{y\in E}A_y\big) f_0\qquad\text{and}\qquad 
\varphi_y=\big(\prod_{\substack{z\in E\\z\neq y}}A_z\big) f_0\quad (y\in E).
\]
By commutativity, $\psi=A_y\varphi_y$ for every $y\in E$.
The functions $\varphi_y$ are non-negative, continuous, compactly supported, and 
satisfy $\sum_{a\in\mathbb Z^n} \varphi_y(a)=1$.
Note that $\psi\in C^1(\mathbb R^n)$ because $\frac{\partial\psi}{\partial x_i}(x)=\varphi_{e_i}(x)-\varphi_{e_i}(x-e_i)$. For $\theta=(\theta_1,\dotsc,\theta_n)\in\{-1,1\}^n$, we thus have
\[
\varphi_\theta (x)-\varphi_\theta(x-\theta)=D_\theta\psi(x)=\sum_{i=1}^n\theta_i D_{e_i}\psi(x)=\sum_{i=1}^n \theta_i \big(\varphi_{e_i}(x)-\varphi_{e_i}(x-e_i)\big).
\]

Let 
\[
M_n=\{0\}\cup \bigcup_{y\in E} \{a,a+y : a\in\mathbb Z^n,\ \varphi_y(a)\neq 0\}.
\]
Since $E$ is finite and all the functions $\varphi_y$ have compact support, $M_n$ is finite. We equip $M_n$ with the metric inherited from $\ell_\infty^n$.
For $i=1,\dotsc,n$, define
\[
u_i=\sum_{a\in\mathbb Z^n}\varphi_{e_i}(a)(\delta_{a+e_i}-\delta_a).
\]
Then 
\[
    \langle u_i,\pi_j\rangle=\sum_{a\in\mathbb Z^n}\varphi_{e_i}(a)\big(\pi_j(a+e_i)-\pi_j(a)\big)=\delta_{ij}\sum_{a\in\mathbb Z^n}\varphi_{e_i}(a)=\delta_{ij}.
\]
For $\theta=(\theta_1,\dotsc,\theta_n)\in\{-1,1\}^n$, 
\begin{align*}
    \sum_{i=1}^n\theta_i u_i
    &=\sum_{a\in\mathbb Z^n}\sum_{i=1}^n\theta_i\big(\varphi_{e_i}(a-e_i)-\varphi_{e_i}(a)\big)\delta_a\\
    &=\sum_{a\in\mathbb Z^n} \big(\varphi_\theta(a-\theta)-\varphi_\theta(a)\big)\delta_a\\
    &=\sum_{a\in\mathbb Z^n}\varphi_\theta(a)(\delta_{a+\theta}-\delta_a).
\end{align*}
Hence
\[
\big\|\sum_{i=1}^n\theta_i u_i\big\|\leq\sum_{a\in\mathbb Z^n}\varphi_\theta(a)\,\|\delta_{a+\theta}-\delta_a\|= 1.
\]
Thus, the two sufficient conditions are satisfied, and the theorem follows.

\section{Further remarks}
\subsection*{1. Norm-one projections}
The construction yields explicit norm-one projections on both $\mathcal F(M_n)$ and its dual. Put $Y=\operatorname{span}\{u_1,\dotsc,u_n\}$ and define $P\mu=\sum_{j=1}^n\langle\mu,\pi_j\rangle u_j$, $\mu\in\mathcal F(M_n)$. Then $Pu_i=u_i$ for every $i=1,\dotsc,n$. And by~\eqref{eq:1},
\[
\|P\mu\|=\max_{1\leq j\leq n}|\langle\mu,\pi_j\rangle|\leq \|\mu\|.
\]
Hence, $P$ is a norm-one projection from $\mathcal F(M_n)$ onto $Y$. 

There is a corresponding dual statement. Let $Z=\operatorname{span}\{\pi_1,\dotsc,\pi_n\}\subset \operatorname{Lip}_0(M_n)$. We have $\|\sum_{j=1}^n a_j\pi_j\|_{\operatorname{Lip}}\leq \sum_{j=1}^n |a_j|$ for all $a_1,\dotsc,a_n\in\mathbb R$. On the other hand, by \eqref{eq:1} and $\langle u_i,\pi_j\rangle=\delta_{ij}$, 
\begin{align*}
    \big\|\sum_{j=1}^n a_j\pi_j\big\|_{\operatorname{Lip}}&\geq \max_{\theta\in\{-1,1\}^n}\big|\big\langle \sum_{i=1}^n\theta_i u_i,\sum_{j=1}^n a_j\pi_j\big\rangle\big|=\max_{\theta\in\{-1,1\}^n}\big|\sum_{j=1}^n\theta_ja_j\big|=\sum_{j=1}^n|a_j|.
\end{align*}
Thus, $Z$ is linearly isometric to $\ell_1^n$. 

Define $Qf=\sum_{j=1}^n\langle u_j,f\rangle\pi_j$, $f\in \operatorname{Lip}_0(M_n)$. Then $Q\pi_j=\pi_j$ for every $j=1,\dotsc,n$, and 
\begin{align*}
    \|Qf\|_{\operatorname{Lip}}&=\sum_{i=1}^n|\langle u_i,f\rangle|=\max_{\theta\in\{-1,1\}^n}\big|\big\langle\sum_{i=1}^n \theta_i u_i,f\big\rangle\big|\leq \|f\|_{\operatorname{Lip}}.
\end{align*}
Hence, $Q$ is a norm-one projection from $\operatorname{Lip}_0(M_n)$ onto $Z$. In particular, $\operatorname{Lip}_0(M_n)$ contains a $1$-complemented linear isometric copy of $\ell_1^n$. 

\subsection*{2. A universal countable compact metric space.} 
There exists a single countable compact metric space whose Lipschitz-free space contains a $1$-complemented copy of $\ell_\infty^n$ for every $n$. For each $n$, choose 
$\lambda_n>0$ such that $\lambda_n\max_{x\in M_n}d(x,0)\leq 2^{-n}$, and let $K_n$ be a  copy of $M_n$ with its metric multiplied by $\lambda_n$. Form the pointed metric wedge 
$K=\bigvee_{n=1}^\infty K_n$. Thus, the base points of the spaces $K_n$ are identified, distances within each $K_n$ remain unchanged, and $d(x,y)=d(x,0)+d(y,0)$ whenever $x\in K_n$ and $y\in K_m$ with $n\neq m$.

The space $K$ is countable and compact. 
For every $n$, the map $R_n\colon K\to K_n$ given by 
\[R_n(x)= \begin{cases}
x, & x\in K_n,\\
0, & x\notin K_n,
\end{cases}
\] is a $1$-Lipschitz retraction. Its linearisation is therefore a norm-one projection
from $\mathcal F(K)$ onto $\mathcal F(K_n)$. Scaling the metric does not change the linear isometry class of the corresponding Lipschitz-free space. Hence, $\mathcal F(K)$ contains a $1$-complemented linear isometric copy of $\ell_\infty^n$ for every
$n\in\mathbb N$. It follows that $\mathcal F(K)$ contains a linear isometric copy of
every finite-dimensional polyhedral real Banach space.

\subsection*{3. Minimal cardinality and simple constructions} 
For $n\in\mathbb N$, let 
\[m(n)=\min\{|M|:\mathcal F(M)\text{ contains a linear isometric copy of }\ell_\infty^n\}.\] 
We have $m(1)=2$ and $m(2)=3$. In the construction above, it is enough to use one representative from each antipodal pair in $\{-1,1\}^n$ because the estimate for $-\theta$ is the same as that for $\theta$. As $\operatorname{supp} f_0\subset [-1,1]^n$ and $A_y$ adds at most the segment $[0,y]$ to the support, this gives $m(n)\leq (2^{n-1}+4)^n$.

A disjoint-roadmap argument as in \cite[Proposition~3.3]{MR4455126} gives $m(n)\geq 2^{n-1}$. We omit the details.
Corollary~3.8 of \cite{MR4181163} gives a four-point metric space $M$ such that $\mathcal F(M)$ is isometric to $\ell_\infty^3$, while the eight-point example in \cite[Section~8]{MR4301142} gives an eight-point metric space whose Lipschitz-free space contains an isometric copy of $\ell_\infty^4$. Hence $m(3)=4$ and $m(4)=8$. Determining the growth of $m(n)$ and finding simpler constructions for $n\geq 5$ remain open.

\section*{Acknowledgements}
The author thanks Mikhail Ostrovskii for his continued interest in this work and for his helpful and encouraging comments. The author also thanks Andre Ostrak for valuable discussions during the development of this work, and Tomáš Raunig for showing him an AI-generated example for $n=7$, which helped motivate the search for a general construction.

\bibliographystyle{amsplain}
\bibliography{references}
\end{document}